\documentclass[a4paper,12pt]{article}

\usepackage{amsmath}                                                            
\usepackage{amssymb}                                                            
\usepackage{amsfonts}

\usepackage{mathrsfs}

\newcommand{\ZZ}{\mathbb{Z}}

\DeclareFontFamily{OT1}{pzc}{}
\DeclareFontShape{OT1}{pzc}{m}{it}{<-> s * [1.200] pzcmi7t}{}
\DeclareMathAlphabet{\mathpzc}{OT1}{pzc}{m}{it}

 \newcommand{\Icg}[2]{\mathrm{ICG}({#2},{#1})}   
\newcommand{\Ene}[2]{\mathpzc{E}({#2},{#1})}   

\newcommand{\Emin}[1]{\mathpzc{E}_{\mathrm{min}}({#1})}   
\newcommand{\Emax}[1]{\mathpzc{E}_{\mathrm{max}}({#1})}   
\newcommand{\Emaxstrich}[1]{{\mathpzc{E}}\sp{\ast}_{\mathrm{max}}({#1})}   
\newcommand{\Emaxr}[2]{\mathpzc{E}_{\mathrm{max}}({#1},{#2})}

\def \C{\hbox{\sf \rlap{\kern.25em \vrule width.1em height1.6ex depth-.1ex}C}}
\def \D{{\sf I\kern-1.5ptD\,}}
\def \K{{\sf I\kern-1.5ptK\,}}
\def \N{{\sf I\kern-1.5ptN\,}}
\def \P{{\sf I\kern-1.5ptP\,}}
\def \Q{\hbox{\sf \rlap{\kern.25em \vrule width.1em height1.6ex depth-.1ex}Q}}
\def \R{{\sf I\kern-1.5ptR\,}}
\def \T{{\sf T\kern-6.5ptT\,}}
\def \Z{{\sf Z\kern-5.0ptZ\,}}

\def\begitem#1 {\bigskip\pagebreak[1]%
     \refstepcounter{subsection}{\nopagebreak[4]%
     \thesubsection\hskip 0.5truecm}
     {\sc#1}\hskip 1pt.\nopagebreak[4]\par\nopagebreak[4]%
      \begin{enumerate}\rm\nopagebreak[4]}
\def\BEGITEMK#1 #2{\bigskip\pagebreak[1]%
      \refstepcounter{subsection}{\nopagebreak[4]%
     \thesubsection\hskip 0.5truecm}\nopagebreak[4]
     {\bf#1}\hskip 1pt.\nopagebreak[4]\par\nopagebreak[4]%
     \medskip\nopagebreak[4]\rm#2\nopagebreak[4]%
     \begin{enumerate}\nopagebreak[4]\rm}
\def\enditem{\end{enumerate}}
\newtheorem{theorem}{Theorem}[section]
\newtheorem{corollary}{Corollary}[section]

\newtheorem{proposition}{Proposition}[section]

\begin{document}

\title{Integral circulant graphs of prime power order with maximal energy}
\author{J.W.~Sander und T.~Sander}
\maketitle

\begin{abstract}

The energy of a graph is the sum of the moduli of the eigenvalues of its adjacency matrix. We study the energy of integral circulant graphs, also called gcd graphs,
which can be characterized by their vertex count $n$ and a set $\cal D$ of divisors of $n$ in such a way that they have vertex set $\mathbb{Z}_n$ and edge set $\{\{a,b\}:\, a,b\in\mathbb{Z}_n,\, \gcd(a-b,n)\in {\cal D}\}$.
Using tools from convex optimization, we study the maximal energy among all integral circulant graphs of prime power order $p^s$ and varying divisor sets $\cal D$. 
Our main result states that this maximal energy approximately lies between $s(p-1)p^{s-1}$ and twice this value. We construct
suitable divisor sets for which the energy lies in this interval.
We also characterize hyperenergetic integral circulant graphs of prime power order and exhibit an interesting topological property of their divisor sets.
\end{abstract}

\\

\\

\section{Introduction}

Concerning the energies of integral circulant graphs, an interesting open problem is the characterization of those graphs having 
maximal energy among all integral circulant graphs with the same given number of vertices. 
The goal of this paper is to establish clarity concerning this question, for integral circulant graphs of prime power order,  by showing how to construct 
 such graphs with a prescribed number of vertices whose energy comes close to the desired maximum.
In the course of this, we approximately determine the maximal energy itself. We rely on 
a closed formula for the energy of an integral circulant graph with prime power order that was established in \cite{SA2}.

A circulant graph is a graph whose adjacency matrix (with respect to a suitable vertex indexing) can be constructed from its first row by a process of continued rotation of entries. An integral circulant graph is a circulant graph whose adjacency matrix has only integer eigenvalues.
The integral circulant graphs belong to the class of Cayley graphs. By a result of {\sc So} \cite{SO}, they are in fact exactly the class of the so-called gcd graphs \cite{SO}, a class that originally arose as a generalization of unitary Cayley graphs. 
The gcd graphs have first been described by {\sc Klotz} and {\sc Sander} in \cite{KLO} and further studied e.g.~by
{\sc Ba\v si\'c} and {\sc Ili\'c} \cite{BAS}, \cite{ILI}. The way the gcd graphs are defined serves us well, so throughout this paper we shall make use of this particular perspective of perceiving integral circulant graphs.
Given an integer $n$ and a set ${\cal D}$ of positive divisors of $n$, the integral circulant graph $\Icg{\cal D}{n}$ is  
defined as the corresponding gcd graph having vertex set $\ZZ_n=\{0,1,\ldots,n-1\}$ and edge set $\{ \{a,b\}:~a,b\in \ZZ_n,~ \gcd(a-b,n)\in {\cal D}\}$.
 We consider only loopfree gcd graphs, i.e. $n\notin {\cal D}$.

The \textit{energy} $E(G)$ of a graph $G$ on $n$ vertices is defined as
\[
     E(G)=\sum_{i=1}^n \vert\lambda_i\vert,
\]
 where $\lambda_1,\ldots,\lambda_n$ are the
eigenvalues of the adjacency matrix of $G$. This concept has been introduced several decades ago
by {\sc Gutman} \cite{GUT}, and with 
slight modification it can even be extended to arbitrary real rectangular matrices, cf.~\cite{NIK} and \cite{KHA}.

There exist many bounds for the graph energy, see {\sc Brualdi} \cite{BRU} for a 
short survey. One example is the bound
\[
 E(G) \leq \frac{n}{2} \left( \sqrt{n} + 1\right)
\]
due to {\sc Koolen} and {\sc Moulton} \cite{KOO} for any graph with $n$ vertices.
There exist infinitely many graphs that achieve this bound. 
If we consider only the class of circulant graphs, then the question arises how
close one can get to this bound. {\sc Shparlinksi} \cite{SHP} has given an explicit construction
of an infinite family of graphs that asymptotically achieves the bound.

Another well-known result is due to {\sc Balakrishnan}, who gives an upper bound $B=k+\sqrt{k(n-1)(n-k)}$ for the energy of
a $k$-regular graph on $n$ vertices (see \cite{BAL}). 
{\sc Li et al.} \cite{LI} have shown that for every $\varepsilon > 0$
one can actually find infinitely many $k$-regular graphs $G$ such
that $\frac{E(G)}{B}>1-\varepsilon$.

There has been some recent work on the energy of unitary Cayley graphs, which are exactly the gcd graphs with ${\cal D}=\{1\}$.
Let us abbreviate $\Ene{\cal D}{n}=E(\Icg{\cal D}{n})$ and let $n=p_1^{s_1}\cdots p_k^{s_k}$.
Then, in the context of gcd graphs, the following result has been obtained
by {\sc Ramaswamy} and {\sc Veena} \cite{RAM} and, independently, by {\sc Ili\'c} \cite{ILI}:
\[
 \Ene{\{1\}}{n} = 2^{k}\varphi(n),
\]
where $\varphi$ denotes Euler's totient function.

{\sc Ili\'c} \cite{ILI} has slightly generalized this to some gcd graphs that are not
unitary Cayley graphs:
\[
\begin{array}{lll}
  \Ene{\{1,p_i\}}{n}   &= 2^{k-1} p_i \varphi(n/p_i), &\qquad \mbox{provided that}~s_i=1, \\
  \Ene{\{p_i,p_j\}}{n} &= 2^k \varphi(n),              &\qquad \mbox{provided that}~s_1=\ldots=s_k=1.  
\end{array}
\]

In \cite{SA2} the authors proved an explicit formula for the energy of $\Icg{\cal D}{p^s}$
for any prime power $p^s$ and any divisor set ${\cal D} = \{ p^{a_1}, p^{a_2} , \ldots , p^{a_r}\}$ with
$0\le a_1<a_2< \ldots < a_r \le s-1$, namely
\begin{equation}
\Ene{{\cal D}}{p^s} = 2(p-1)\left(p^{s-1}r -(p-1) \sum_{k=1}^{r-1} \sum_{i=k+1}^r p^{s-a_i+a_k-1}\right). \label{explform}
\end{equation}

The study of energies is usually linked to the search for hyperenergetic graphs.
A graph $G$ on $n$ vertices is called \textit{hyperenergetic} if its energy is greater than the energy
of the complete graph on the same number of vertices, i.e.~if 
$E(G) > E(K_n) = 2(n-1)$. Initially, the existence of hyperenergetic graphs had been doubted, but
then more and more classes of hyperenergetic graphs were discovered. For example, {\sc Hou} and {\sc Gutman} show in \cite{HOU} that if a 
graph $G$ has more than $2n-1$ edges, then its line graph $L(G)$ is necessarily hyperenergetic. 
Consequently, $L(K_n)$ is 
hyperenergetic for all $n\geq 5$, a fact that seems to have been known before.

Work by {\sc Stevanovi\'c} and {\sc Stankovi\'c} \cite{STE} indicates that 
the class of circulant graphs contains a wealth of hyperenergetic graphs.
Although integral circulant graphs are quite rare among circulant graphs (cf. \cite{AHM}),
the subclass of integral circulant graphs still exhibits many hyperenergetic members.
For example, it has been shown by {\sc Ramaswamy} and {\sc Veena} \cite{RAM}
that almost all unitary Cayley graphs
on $n$ vertices are hyperenergetic. The necessary and sufficient 
condition is that  $n$ has at least $3$
distinct prime divisors or that $n$ is odd in case of only two prime divisors.
Consequently, there exist no gcd graphs $\Icg{\cal D}{p^s}$ with ${\cal D}=\{1\}$  that are hyperenergetic. 
However, for less trivial divisor sets it is also possible to find 
hyperenergetic gcd graphs on $p^s$ vertices. Some examples are given in \cite{SA2}. For instance,
for $p\geq 3$ and $s\geq 3$, the choice ${\cal D}=\{1, p^{s-1}\}$ yields a hyperenergetic gcd graph.

Not surprisingly, the class of graphs 
$\Icg{\cal D}{p^s}$  contains also non-hyperenergetic elements, termed \textit{hypoenergetic}.
For the minimal energy $\Emin{n}$ of all integral circulant graphs with $n$ vertices it has been shown in \cite{SA2} that
$$\Emin{p^s} = 2(p-1)p^{s-1} = E(K_{p^s})-E(K_{p^{s-1}}).$$
This follows directly from equation (\ref{explform}). 
The minimal energy is achieved exactly for the singleton divisor sets. 

The maximum energy of graphs $\Icg{\cal D}{p^s}$ is not as easily described. A classification 
of integral circulant graphs of prime power order $p^s$ with very small exponent having maximal energy has been provided in \cite{SA2}, 
but it became clear that a general result as simple as in the case of minimal
energy could not be expected. It will be the purpose of this article to clarify the structure of divisor sets
imposing maximal energy on the corresponding gcd graph. 
Our main Theorem \ref{t5} states that the maximal energy 
among all integral circulant graphs of prime power order $p^s$ and varying divisor sets $\cal D$
approximately lies between $s(p-1)p^{s-1}$ and twice this value.
Tools from convex optimization will turn out to be the
appropriate machinery to reach that goal. We shall compute bounds for the maximum energy
and describe how to construct divisor sets for integral circulant graphs on $p^s$ vertices that have near
maximal energy.
Along the way, we characterize hyperenergetic integral circulant graphs of prime power order and 
exhibit an interesting topological property of their divisor sets. Namely, the set containing all ordered exponent tuples 
corresponding to these divisor sets can be obtained by intersecting an integer grid with a suitable convex set.

%
%

\section{Preliminary definitions and results}

For any positive integer $n$, let
$$ \Emax{n} := \max\,\{ \Ene{\cal D}{n}:\;\; {\cal D}\subseteq \{1\le d<n:\; d\mid n\} \}.$$

For given ${\cal D} = \{ p^{a_1}, p^{a_2} , \ldots , p^{a_r}\}$ with
$0\le a_1 < \ldots < a_r \le s-1$, we have by (\ref{explform}), i.e. by Theorem 2.1 in \cite{SA2}, that
\begin{equation}
\Ene{{\cal D}}{p^s} = 2(p-1)p^{s-1}\left(r -(p-1)h_p(a_1,\ldots,a_r) \right),  \label{ft3}
\end{equation}
where
$$ h_p(x_1,\ldots,x_r) := \sum_{k=1}^{r-1}\sum_{i=k+1}^r \,\frac{1}{p^{x_i-x_k}}$$
for arbitrary real numbers $x_1, \ldots, x_r$.
In order to evaluate $\Emax{p^s}$, our main task will be to determine 
$$ \Emaxr{p^s}{r}:= \max\, \{\Ene{{\cal D}}{p^s}:\; {\cal D}\subseteq \{1\le d<n:\; d\mid n\}, \; \vert {\cal D} \vert =r \} $$
as precisely as possible, given a fixed integer $r$. Therefore, we define for $1 \le r \le s+1$ 
\begin{equation}
m_p(s,r) := \min\,\{ h_p(a_1,\ldots,a_r): \; 0\le a_1 < a_2 < \ldots < a_r\le s \mbox{ with $a_i\in\mathbb{Z}$} \}. \label{intmin}
\end{equation}
It is then clear from (\ref{ft3}) that
\begin{equation}
\Emaxr{p^s}{r} = 2(p-1)p^{s-1}\left(r -(p-1)m_p(s-1,r)\right).   \label{ft4}
\end{equation}
Later on it remains to compute
\begin{equation}
\Emax{p^s} = \max\, \{ \Emaxr{p^s}{r}: \; 1\le r \le s \}. \label{ft5}
\end{equation}

\begin{proposition}{\label{p1}} 
Let $p$ be a prime. Then
\begin{itemize}
\item[(i)] $m_p(s,2) = \frac{1}{p^{s}}$ for all integers $s\ge 1$, and the minimum is attained only for $a_1=0$ and $a_2=s$.
\item[(ii)] $m_p(s,3) = \frac{1}{p^{[s/2]}} + \frac{1}{p^{s}} + \frac{1}{p^{s-[s/2]}}$ for all integers $s\ge 2$.
The minimum is only obtained for $a_1=0$, $a_2=[s/2]$ (or, additionally, for $a_2=[s/2]+1$ if $s$ is odd) and $a_3=s$.
\end{itemize}
\end{proposition}

{\sc Proof. } Proposition 3.1 in \cite{SA2}.
\vspace{-.6truecm}
\begin{flushright}$\Box$\end{flushright}

A set ${\cal D}\subseteq \{1\le d<n:\; d\mid n\}$ is called \textit{$n$-maximal} if $\Ene{\cal D}{n} =  \Emax{n}$.
As a consequence of Proposition \ref{p1} and some other results in \cite{SA2}, we obtained
\begin{theorem}{\label{t3}} Let $p$ be a prime.
Then 
\begin{itemize}
\item[(i)] $\Emax{p} = 2(p-1)$ with the only $p$-maximal set ${\cal D} = \{1\}$.
\item[(ii)] $\Emax{p^2} = 2(p-1)(p+1)$ with the only $p^2$-maximal set ${\cal D} = \{1,p\}$. 
\item[(iii)] $\Emax{p^3} = 2(p-1)(2p^2-p+1)$ with the only $p^3$-maximal set ${\cal D} = \{1,p^2\}$, except 
for the prime $p=2$ for which ${\cal D} = \{1,2,4\}$ is also $2^3$-maximal.
\item[(iv)] $\Emax{p^4} = 2(p-1)(2p^3+1)$ with the only $p^4$-maximal sets ${\cal D} = \{1,p,p^3\}$ and
${\cal D} = \{1,p^2,p^3\}$.
\end{itemize}
\end{theorem}

{\sc Proof. } Theorem 3.2 in \cite{SA2}.
\vspace{-.6truecm}
\begin{flushright}$\Box$\end{flushright}

One can prove formulae for $\Emax{p^s}$ with arbitrary 
exponent $s$ by using (\ref{ft4}) and (\ref{ft5}).
As indicated in Proposition \ref{p1}, we need to choose integers
$0\le a_1 \le a_2 \le \ldots \le a_r \le s-1$ in such a way that
$$ h_p(a_1,\ldots,a_r) = \sum_{k=1}^{r-1}\sum_{i=k+1}^r \,\frac{1}{p^{a_i-a_k}}$$
becomes minimal. The choice of $a_1=0$ and $a_r=s-1$ is clearly compulsory. 

The case $r=3$ (cf. Prop. \ref{p1}(ii)) suggests to place $a_1,a_2,\ldots,a_{r-1},a_r$ equidistant in the interval $[0,s-1]$.
A minor obstacle is the fact that 
the corresponding choice $a_i:=\frac{(i-1)(s-1)}{r-1}$ ($1\le i \le r$) does not yield integral numbers as required. Taking 
nearest integers easily resolves this problem, but only at the cost of approximate instead of exact formulae.
More seriously, it turns out that in general, even allowing real $a_i$, their 
equidistant positioning does not yield the desired minimum $m_p(s,r)$. 
The cases presented in Proposition \ref{p1} do not yet exhibit this problem since it
makes its debut for $r=4$. An illuminating example can be found in the final section of \cite{SA2}.

For the sake of being able to use analytic methods, we define for a prime $p$, a positive real number $\sigma$ and a positive integer $r$
\begin{equation} 
\tilde{m}_p(\sigma,r) := \min\,\{ h_p(\alpha_1,\ldots,\alpha_r): \; 0\le \alpha_1 \le \alpha_2\le \ldots \le \alpha_r \le \sigma, \; \alpha_i \in \mathbb{R} \}. \label{realmin}
\end{equation}
Observe that now the $\alpha_i$ may be real numbers as opposed to integers in the definition of $m_p(s,r)$. It is obvious that for $r\ge 2$
\begin{equation}
 \tilde{m}_p(\sigma,r) = \min\,\{ h_p(0,\alpha_2,\ldots,\alpha_{r-1},\sigma): \; 0\le \alpha_2 \le \ldots \le \alpha_{r-1} \le \sigma, \; \alpha_i \in \mathbb{R} \}.
\label{mtilde}
\end{equation}
Clearly, $\tilde{m}_p(\sigma,2)=1/p^{\sigma}$, uniquely obtained for $\alpha_1=0$, $\alpha_2=\sigma$, and  
$\tilde{m}_p(\sigma,3)=1/p^{\sigma}+2/p^{\sigma/2}$, uniquely obtained for $\alpha_1=0$, $\alpha_2=\sigma/2$, $\alpha_3=\sigma$ (cf. Proposition \ref{p1}(ii)).

%
%

\section{Tools from convex optimization}

In order to determine $\tilde{m}_p(\sigma,r)$ in general it is crucial to observe that $h_p(x_1,\ldots,x_r)$ is a convex function.
\begin{proposition}{\label{p2convex}} 
Let $r$ be a fixed positive integer, $b\neq 1$ a fixed positive real number and $p$ a fixed prime.
\vspace{-.07in}
\begin{itemize}
\item[(i)] The real function
$$  g_b(y_1,\ldots,y_r) :=  \sum_{k=1}^r \sum_{i=k}^r \prod_{j=k}^i \frac{1}{b^{y_j}}$$ 
is strictly convex on $\mathbb{R}^r$.
\item[(ii)] The function $h_p(x_1,\ldots,x_r)$ is convex on $\mathbb{R}^r$.
\end{itemize}
\end{proposition}

{\sc Proof. }\newline 
(i)  Let $(u_1,\ldots,u_r)\neq (v_1,\ldots,v_r)$ be arbitrary real vectors.
On setting
$$ U_{k,i}:= \prod_{j=k}^i \frac{1}{b^{u_j}} \quad \mbox{ and } \quad V_{k,i}:= \prod_{j=k}^i \frac{1}{b^{v_j}},$$
we have $U_{k,k}=b^{-u_k}$ and $V_{k,k}=b^{-v_k}$ for $1\le k \le r$. Since $0< b\neq 1$ and $u_k\neq v_k$ for at least
one $k$, we have $U_{k,k}\neq V_{k,k}$ for that $k$. 
By the inequality between the weighted arithmetic and the weighted geometric mean, which is an immediate consequence of Jensen's inequality (cf. \cite{ROB1}, p. 1100, Thms. 17 and 18), we have $U^t \cdot V^{1-t} \le tU+(1-t)V$ for
all positive real numbers $U$ and $V$ and all $0<t<1$, with equality if and only if $U=V$.
It follows that for all $1\le k \le i \le r$
$$ U_{k,i}^t V_{k,i}^{1-t} \le t U_{k,i} + (1-t) V_{k,i},$$
and strict inequality for at least one pair $k,i$.
Hence
\begin{eqnarray*}
g_b\big(t(u_1,\ldots,u_r)+(1-t)(v_1,\ldots,v_r)\big) &=& \sum_{k=1}^{r} \sum_{i=k}^r  U_{k,i}^t V_{k,i}^{1-t} \\
&<& \sum_{k=1}^{r} \sum_{i=k}^r t  U_{k,i} + \sum_{k=1}^{r} \sum_{i=k}^r (1-t) V_{k,i}\\
&=&  t g_b(u_1,\ldots,u_r)+(1-t)g_b(v_1,\ldots,v_r).
\end{eqnarray*}
This proves that $g$ is strictly convex on $\mathbb{R}^r$.
\vspace{.2in}\newline
(ii) Let $(u_1,\ldots,u_r), (v_1,\ldots,v_r)\in \mathbb{R}^r$ and $0<t<1$. Then by (i)
\begin{eqnarray*}
& & h_p\big(t(u_1,\ldots,u_r)+(1-t)(v_1,\ldots,v_r)\big) = \\
&=& g_p\big(t(u_2-u_1,u_3-u_2,\ldots,u_r-u_{r-1})+(1-t)(v_2-v_1,v_3-v_2,\ldots,v_r-v_{r-1})\big) \\
&\le& t g_p\big(u_2-u_1,u_3-u_2,\ldots,u_r-u_{r-1}\big)+(1-t)g_p\big(v_2-v_1,v_3-v_2,\ldots,v_r-v_{r-1}\big)\\
&=& t h_p(u_1,\ldots,u_r)+(1-t)h_p(v_1,\ldots,v_r),
\end{eqnarray*}
which shows the convexity of $h_p$.
\vspace{-.6truecm}
\begin{flushright}$\Box$\end{flushright}

An easy corollary of (\ref{explform}) is a characterisation of the hyperenergetic gcd graphs of prime power order, namely that
$\Icg{\cal D}{p^s}$ with
${\cal D} = \{ p^{a_1}, p^{a_2} , \ldots , p^{a_r}\}$ and $0\le a_1<a_2< \ldots < a_r \le s-1$ is hyperenergetic if and only if
\begin{equation}
\sum_{k=1}^{r-1} \sum_{i=k+1}^r \frac{1}{p^{a_i-a_k}} < \frac{1}{p-1} \left(r - \frac{p^s-1}{p^{s-1}(p-1)} \right)  \label{cor2.2}
\end{equation}
(cf. Corollary 2.2 in \cite{SA2}). As a first consequence of the convexity of $h_p$ we are able to refine this by showing that the set of hyperenergetic integral circulant graphs 
has a nice topological feature. Given a prime $p$ and positive integers $r\le s$, we define 
${\cal H}(p^s,r)$ as the set containing all $(a_1,\ldots,a_r)\in \mathbb{Z}^r$ with $0\le a_1<\ldots<a_r\le s-1$ and the property that $\Icg{\{p^{a_1},\ldots,p^{a_r}\}}{p^s}$ is hyperenergetic. Then we can derive the following remarkable statement:
\begin{corollary}{\label{c2a}} Let $p$ be a prime and $r\le s$ positive integers. Then there is a convex set $C\subseteq \mathbb{R}^r$
such that ${\cal H}(p^s,r) = C \cap \mathbb{Z}^r$. 
\end{corollary}

{\sc Proof. }
For fixed $p$, $s$ and $r$, we define
$$ c(p,s,r):= \frac{1}{p-1} \left(r - \frac{p^s-1}{p^{s-1}(p-1)} \right).$$
Since $h_p$ is convex on $\mathbb{R}^r$ by
Proposition \ref{p2convex}, the so-called level set
$$ L:= \{(x_1,\ldots,x_r)\in \mathbb{R}^r: \; h_p(x_1,\ldots,x_r)< c(p,s,r)\}$$
is also convex (cf. \cite{ROC}, p.~8 and Prop. 2.7). Since 
$$ K:= \{(x_1,\ldots,x_r)\in \mathbb{R}^r: \; 0\le x_1 < x_2< \ldots <x_r\le s-1 \}$$
is obviously convex as well, the intersection $C:=L\cap K$ has the same property.
By (\ref{cor2.2}) we know that some $(a_1,\ldots,a_r)\in \mathbb{Z}^r$ lies in ${\cal H}(p^s,r)$ if
and only if $0\le a_1<\ldots<a_r\le s-1$ and $h_p(a_1,\ldots,a_r)< c(p,s,r)$, hence ${\cal H}(p^s,r) = C \cap \mathbb{Z}^r$.
\vspace{-.6truecm}
\begin{flushright}$\Box$\end{flushright}

We shall use some further standard results from convex optimization.
\begin{proposition}{\label{p2}} 
Let $f: U \rightarrow \mathbb{R}$ be a strictly convex function defined on a convex set $U\subseteq \mathbb{R}^r$.
\begin{itemize}
\item[(i)] If $U$ is an open set then each extremal point of $f$ is a minimum.
\item[(ii)] If $f$ has a minimal point on $U$ then it is unique.
\end{itemize}
\end{proposition}

{\sc Proof. }
The proofs of the assertions can be found in \cite{ROB}, pp. 123-124, Theorems A and C,  in \cite{KAP}, or in \cite{ROC}, Thm. 2.6.
\vspace{-.6truecm}
\begin{flushright}$\Box$\end{flushright}

Our main tool for the computation of $\tilde{m}_p(\sigma,r)$ will be 
\begin{proposition}{\label{p3}} 
Let $r\ge 1$ be a fixed integer. We define the real function 
$$ f(x_1,\ldots,x_r) := \sum_{k=1}^{r} \sum_{i=k}^r \prod_{j=k}^i\, x_j$$
for $(x_1,\ldots,x_r)\in [0,1]^r$. Let $0<\rho\le 2^{-r}$. Then
$$ \min\, \{f(x_1,\ldots,x_r): \;    (x_1,\ldots,x_r)\in [0,1]^r, \;\, x_1\cdot x_2\cdots x_r=\rho \} = (r + \mu(\rho,r))\cdot\mu(\rho,r),$$
where $\mu(\rho,r):=\nu(\rho,r)/(1-\nu(\rho,r))$ and $x=\nu(\rho,r)$ is the unique real solution of the equation $x^r = \rho(1-x)^2$ on the
interval $[0,1]$. The minimum obviously equals $\rho$ for $r=1$, and it is $\rho + 2\sqrt{\rho}$ for $r=2$.

There is a unique minimizer for each $r$, namely
\begin{equation*}
\left\{ \begin{array}{ll}
\rho\in [0,1]   &  \mbox{ for $r=1$,} \\
(\mu(\rho,2),\mu(\rho,2)) \in [0,1]^2 & \mbox{ for $r=2$,} \\
(\mu(\rho,r),\nu(\rho,r),\ldots,\nu(\rho,r),\mu(\rho,r)) \in [0,1]^r   & \mbox{ for $r\ge 3$.} 
\end{array} \right.
\end{equation*}
In the special case $r=2$, we have explicitly $\mu(\rho,2)=\sqrt{\rho}$.
\end{proposition}

{\sc Proof. } 
For $r\le 2$, we have to deal with nothing more than quadratic equations, and in these cases all assertions follow easily
from standard analysis.

For $r\ge 3$, we use the method of Lagrange multipliers to obtain necessary conditions for local minima of $f(x_1,\ldots,x_r)$ subject to 
the constraint $x_1\cdot x_2\cdots x_r=\rho$. Accordingly, let
$$ F(x_1,\ldots,x_r,\lambda) := f(x_1,\ldots,x_r) + \lambda (\rho - x_1\cdot x_2\cdots x_r).$$
A necessary condition for a local minimum is that all partial derivatives $F_{x_t}:=\frac{\partial F}{\partial x_t}$ $(1\le t \le r)$ as well as 
$F_{\lambda}:=\frac{\partial F}{\partial \lambda}$
vanish at that point. We have for $1\le t \le r$
\begin{eqnarray*}
 f_{x_t}(x_1,\ldots,x_r) &=&    \sum_{k=1}^r \sum_{i=k}^r \frac{\partial}{\partial x_t} \Bigg(\prod_{j=k}^i\, x_j\Bigg)    \\
&=& \sum_{k=1}^{\min\{r,t\}} \sum_{i=\max\{k,t\}}^r \frac{\partial}{\partial x_t} \Bigg(\prod_{j=k}^i\, x_j\Bigg)  = 
\sum_{k=1}^t \sum_{i=t}^r \prod_{\genfrac{}{}{0pt}{1}{j=k}{j\neq t}}^i\, x_j. 
\end{eqnarray*}
Hence
$$ x_t F_{x_t} = x_t f_{x_t}(x_1,\ldots,x_r) - x_t \Bigg( \lambda \prod_{\genfrac{}{}{0pt}{1}{j=1}{j\neq t}}^r\, x_j \Bigg) = 
\sum_{k=1}^t \sum_{i=t}^r \prod_{j=k}^i\, x_j - \lambda \rho,$$
and we want to find all solutions $(x_1,\ldots,x_r)$ of the following system of equations:
\begin{equation}
\left\{\begin{aligned}
\sum_{k=1}^t \sum_{i=t}^r \prod_{j=k}^i\, x_j &= \lambda \rho \quad\quad (1\le t \le r)\, , \\
x_1\cdots x_r &= \rho\, .
\end{aligned}
\right.
\label{p21}
\end{equation}

From now on we consider $r$ and $\rho$ to be fixed and abbreviate $\mu :=\mu(\rho,r)$ and $\nu :=\nu(\rho,r)$.\vspace{.08in}\newline
\underline{Claim}: One solution of (\ref{p21}) is given by $x_1=x_r=\mu$ and $x_2=\ldots = x_{r-1} =\nu$, where
we have $0<\nu<\mu<1$.
\vspace{.08in}

The real function $h(x):= x^r-\rho(1-x)^2$ is strictly increasing on the interval $[0,1]$ with $h(0)=-\rho$ and $h(1)=1$.
Hence $x^r = \rho(1-x)^2$ has a unique solution on $(0,1)$, which is denoted by $\nu$. Since $\nu^r< \rho< 2^{-r}$, we
even know $\nu < 1/2$. This implies $\nu < \mu <1$.

In order to show that $(\mu,\nu,\nu,\ldots,\nu,\mu)\in(0,1)^r$ satisfies (\ref{p21}), we separate terms containing 
$x_1$ or $x_r$ from the others in the double sum of (\ref{p21}) and obtain for $x_1=x_r=\mu$ and $x_2=\ldots =x_{r-1} = \nu$ 
\begin{eqnarray*}
\sum_{k=1}^t \sum_{i=t}^r \prod_{j=k}^i\, x_j &=& x_1\cdots x_r + \sum_{i=t}^{r-1} \prod_{j=1}^i\, x_j
+ \sum_{k=2}^t  \prod_{j=k}^r\, x_j +  \sum_{k=2}^t \sum_{i=t}^{r-1} \prod_{j=k}^i\, x_j \\
&=& \rho + \sum_{i=t}^{r-1} \mu\cdot \nu^{i-1} + \sum_{k=2}^t \nu^{r-k}\cdot \mu + \sum_{k=2}^t \sum_{i=t}^{r-1} \nu^{i-k+1} \\
&=& \rho + \mu\cdot \frac{\nu^{t-1} - \nu^{r-1}}{1-\nu} + \mu\cdot \frac{\nu^{r-t} - \nu^{r-1}}{1-\nu} + \sum_{k=2}^t \frac{\nu^{t-k+1} - \nu^{r-k+1}}{1-\nu} \\
&=& \rho + \mu\cdot \frac{\nu^{t-1} - 2\nu^{r-1} + \nu^{r-t}}{1-\nu} + \frac{1}{1-\nu}\left(\frac{\nu-\nu^t}{1-\nu} -\frac{\nu^{r-t+1}-\nu^r}{1-\nu}\right)\\
&=& \rho + \frac{1}{(1-\nu)^2}\left( \mu (1-\nu)(\nu^{t-1} - 2\nu^{r-1} + \nu^{r-t}) + (\nu-\nu^t -\nu^{r-t+1}+\nu^r)\right).
\end{eqnarray*}
Since $\mu (1-\nu) = \nu$, we conclude for $(x_1,\ldots,x_r)=(\mu,\nu,\ldots,\nu,\mu)$ that
$$ \sum_{k=1}^t \sum_{i=t}^r \prod_{j=k}^i\, x_j = \rho + \frac{\nu\cdot(1-\nu^{r-1})}{(1-\nu)^2} = \frac{\rho}{\nu^{r-1}}\, .$$
Setting $\lambda := 1/\nu^{r-1}$, this reveals that $(\mu,\nu,\ldots,\nu,\mu)$ satisfies all the upper equations in (\ref{p21}).
The observation that
$$ x_1\cdots x_r = \mu^2\cdot \nu^{r-2} = \frac{\nu^r}{(1-\nu)^2} = \rho$$
completes the proof of the claim.

We now want to show that $f(\mu,\nu,\ldots,\nu,\mu)$ is in fact a minimum subject to the constraint $x_1\cdot x_2\cdots x_r=\rho$, and we shall see as well
that $(\mu,\nu,\ldots,\nu,\mu)$ is the unique minimizer. By Proposition \ref{p2convex}(i) the function 
$ g_2(y_1,\ldots,y_r) = f(2^{-y_1},\ldots,2^{-y_r})$ is strictly convex for all $(y_1,\ldots,y_r)\in \mathbb{R}^r$. 
Our claim has shown that $(\mu,\nu,\ldots,\nu,\mu)$ is an extremal point of $f$ on the set
$\{(x_1,\ldots,x_r)\in [0,1]^r : \; x_1 \cdots x_r = \rho\}$. Therefore $(\mu',\nu',\ldots,\nu',\mu')$ with
$$ \mu':= \frac{\log(1/\mu)}{\log 2} \;\; \mbox{ and } \;\; \nu':= \frac{\log(1/\nu)}{\log 2}$$
is an extremal point of $g_2$ on the set
$$U:= \{(y_1,\ldots,y_r)\in \mathbb{R}_{\ge 0}^r : \; 2^{-y_1}\cdots 2^{-y_r} = \rho \}
   =  \{(y_1,\ldots,y_r)\in \mathbb{R}_{\ge 0}^r : y_1+ \ldots +y_r = \sigma \},$$
where 
$$     \sigma :=  \frac{\log(1/\rho)}{\log 2} \ge r.$$
The set $U$ apparently is the simplex with vertices 
$(\sigma,0,\ldots,0), (0,\sigma,0,\ldots,0),\ldots,(0,\ldots,0,\sigma)$, and therefore a convex subset of $\mathbb{R}_{\ge 0}^r$.
It is immediately clear that $(\mu',\nu',\ldots,\nu',\mu')$ does not lie on the boundary of the simplex $U$, in other
words: the point belongs to the set $U^0$ of inner points of $U$. Altogether the function $g_2$ is strictly convex on the open convex set
$U^0$ and $(\mu',\nu',\ldots,\nu',\mu')$ is an extremal point of $g_2$ in $U^0$. By Proposition \ref{p2}(i) this point
$(\mu',\nu',\ldots,\nu',\mu')$ is a minimal point of $g_2$, and by Proposition \ref{p2}(ii) it is unique. Since $\log_2$ is strictly monotonic,
the point $(\mu,\nu,\ldots,\nu,\mu)$ is the unique minimizer with respect to $f$ on 
$ \{(x_1,\ldots,x_r)\in [0,1]^r : x_1\cdots x_r = \rho \}$.

It remains to calculate the minimum. We obtain similarly as before
\begin{eqnarray*}
f(x_1,\ldots,x_r) &=&  x_1\cdots x_r + \sum_{i=1}^{r-1} \prod_{j=1}^i\, x_j
+ \sum_{k=2}^{r}  \prod_{j=k}^r\, x_j +  \sum_{k=2}^{r-1} \sum_{i=k}^{r-1} \prod_{j=k}^i\, x_j .
\end{eqnarray*}
By evaluating the geometric sums and using the identity $\mu=\nu/(1-\nu)$, it follows that
\begin{eqnarray*}
f(\mu,\nu,\ldots,\nu,\mu)
&=& \mu^2\nu^{r-2}+ 2\mu\cdot \frac{1 - \nu^{r-1}}{1-\nu} + \frac{1}{1-\nu}\left( (r-2)\nu - \frac{\nu^2 - \nu^r}{1-\nu}\right) \\
&=& \frac{\nu^r}{(1-\nu)^2} + \frac{2\nu(1 - \nu^{r-1})}{(1-\nu)^2} + \frac{\nu}{1-\nu}\cdot(r-2) - \frac{\nu^2 - \nu^r}{(1-\nu)^2}  \\
&=& \frac{2\nu-\nu^2}{(1-\nu)^2} + \frac{\nu}{1-\nu}\cdot(r-2)\\
&=& \frac{2\nu}{1-\nu} + \frac{\nu^2}{(1-\nu)^2} + \frac{\nu}{1-\nu}\cdot(r-2)\\
&=& 2\mu +\mu^2 +\mu (r-2) = (r+\mu)\cdot\mu.
\end{eqnarray*}
\vspace{-.9truecm}
\begin{flushright}$\Box$\end{flushright}

\begin{corollary}{\label{c3}} 
Let $p$ be a prime and $r\ge 2$ an integer. For a given real number $\mbox{$\sigma\ge (r-1) \log 2/ \log p$}$ let $x=\tilde{\nu}_p(\sigma,r)$ be the
unique real solution of the equation $p^{\sigma}x^{r-1} = (1-x)^2$ on the interval $[0,1]$, and $\tilde{\mu}_p(\sigma,r):= \tilde{\nu}_p(\sigma,r)/(1- \tilde{\nu}_p(\sigma,r))$.
Then
$$   \tilde{m}_p(\sigma,r) = (r-1+\tilde{\mu}_p(\sigma,r))\cdot\tilde{\mu}_p(\sigma,r),$$
and this value is exclusively attained by
$h_p(\alpha_1,\ldots,\alpha_r)$ for $\alpha_j=\alpha_j(\sigma,r)$ ($1\le j\le r$) defined as $\alpha_1(\sigma,r):=0$, $\alpha_r(\sigma,r):=\sigma$ and                
\begin{equation}
\alpha_j(\sigma,r):= \frac{\log\big(\tilde{\mu}_p(\sigma,r)^{-1}\big)}{\log p} + (j-2)\frac{\log\big(\tilde{\nu}_p(\sigma,r)^{-1}\big)}{\log p} \;\;\quad (r\ge 3; \;\; 2\le j\le r-1).  \label{alphas}
\end{equation}
\end{corollary}

{\sc Proof. }
Let $0\le \alpha_1\le \alpha_2 \le \ldots \le \alpha_r \le \sigma$ be arbitrary, and set
$y_j:= \alpha_{j+1}-\alpha_j$ for $1\le j\le r-1$. Hence $\alpha_i - \alpha_k = y_k + y_{k+1} + \ldots + y_{i-1}$ for
$1\le k < i \le r$. This implies
\begin{eqnarray*}
h_p(\alpha_1,\ldots,\alpha_r) &=& \sum_{k=1}^{r-1} \sum_{i=k+1}^r \frac{1}{p^{y_k+\ldots +y_{i-1}}} \\
&=& \sum_{k=1}^{r-1} \sum_{i=k+1}^r \prod_{j=k}^{i-1} \frac{1}{p^{y_j}} = \sum_{k=1}^{r-1} \sum_{i=k}^{r-1} \prod_{j=k}^i \frac{1}{p^{y_j}}.
\end{eqnarray*}
On setting $x_j:=p^{-y_j}$ for $1\le j \le r-1$, we have $h_p(\alpha_1,\ldots,\alpha_r) = f(x_1,\ldots,x_{r-1})$ for the function $f$ as defined in
Proposition \ref{p3}. By hypothesis, $r-1\ge 1$ and $(x_1,\ldots,x_{r-1})\in [0,1]^{r-1}$. Now we search for conditions to be imposed on
the $\alpha_j$ in order to hit the minimum $\tilde{m}_p(\sigma,r)$. First of all,
we necessarily have $\alpha_1=0$ and $\alpha_r=\sigma$ according to (\ref{mtilde}). Hence
$$ x_1\cdot x_2 \cdot \ldots \cdot x_{r-1} = \frac{1}{p^{y_1+\ldots +y_{r-1}}} = \frac{1}{p^{\alpha_r-\alpha_1}} = \frac{1}{p^{\sigma}}.$$
Again by hypothesis
$$ 0 < \rho := \frac{1}{p^{\sigma}} \le \frac{1}{2^{r-1}}.$$
Applying Proposition \ref{p3}, we conclude that 
$$  \tilde{m}_p(\sigma,r) = (r-1 + \mu(\rho,r-1))\cdot \mu(\rho,r-1) = (r-1 + \tilde{\mu}_p(\sigma,r))\cdot \tilde{\mu}_p(\sigma,r),$$
where this minimum is
exclusively obtained for $x_1=x_{r-1} = \tilde{\mu}_p(\sigma,r)$ and $x_2=x_3=\ldots = x_{r-2}= \tilde{\nu}_p(\sigma,r)$.
This yields for the given values of the $x_j$
$$ \frac{1}{p^{\alpha_2}} = \frac{1}{p^{\alpha_2-\alpha_1}} =\frac{1}{p^{y_1}}= x_1 = \tilde{\mu}_p(\sigma,r),$$
hence $\alpha_2=\log(1/\tilde{\mu}_p(\sigma,r))/\log p$, and for $2\le j \le r-2$
$$ \frac{1}{p^{\alpha_{j+1}-\alpha_j}} =\frac{1}{p^{y_j}}= x_j = \tilde{\nu}_p(\sigma,r),$$
which implies (\ref{alphas}).
\vspace{-.6truecm}
\begin{flushright}$\Box$\end{flushright}

%
%

\section{Integral circulant graphs with maximal energy}

Up to this point, all we have done with respect to general integral circulant graphs with maximal energy refers to
\textit{real} parameters $\alpha_j $ in $h_p(\alpha_1,\ldots,\alpha_r)$. As a consequence, we have the
following upper bound for $\Emaxr{p^s}{r}$, but we are left with
the task to find out how close we can get to the ``real maximum" if we restrict ourselves
to \textit{integral} parameters $a_1,\ldots,a_r$, as required by our problem.

\begin{theorem}{\label{t4}} 
For a prime $p$ and integers $2\le r \le s$, we have 
$$\Emaxr{p^s}{r} \le 2(p-1)p^{s-1} \Big(r-(p-1)\big(r-1+\tilde{\mu}_p(s-1,r)\big)\tilde{\mu}_p(s-1,r)\Big),$$
where $\tilde{\mu}_p$ is defined in Corollary \ref{c3}.
\end{theorem}

{\sc Proof. }
By Corollary \ref{c3} and the definitions of $m_p(s-1,r)$ and $\tilde{m}_p(\sigma,r)$, we immediately have
for any integer $s\ge r-1$
$$ m_p(s-1,r) \ge (r-1+ \tilde{\mu}_p(s-1,r)) \cdot \tilde{\mu}_p(s-1,r).$$
Now our theorem follows at once from this and (\ref{ft4}).

\vspace{-.6truecm}
\begin{flushright}$\Box$\end{flushright}

The first step we take towards integrality of the parameters is to approximate the numbers $\tilde{\mu}_p(s-1,r)$,
$\tilde{\nu}_p(s-1,r)$ and the corresponding $\alpha_j(s-1,r)$, all defined in Corollary \ref{c3}, by simpler terms.

\begin{proposition}{\label{p51}} 
For a prime $p$ and integers $3\le r\le s$, let $\delta:=p^{-\frac{s-1}{r-1}}$. Then we have 
\begin{itemize}
\item[(i)] $\delta \le \tilde{\mu}_p(s-1,r) < \delta + \frac{\delta^2}{1-\delta}$;
\item[(ii)] $\delta - \frac{\delta^2}{1+\delta} \le \tilde{\nu}_p(s-1,r) < \delta \le \frac{1}{p}$;
\item[(iii)] $0 < \log\big(\tilde{\nu}_p(s-1,r)^{-1}\big) - \frac{s-1}{r-1}\log p < \frac{3}{(r-1)p}$;
\item[(iv)] $-\frac{3}{2p} < \log\big(\tilde{\mu}_p(s-1,r)^{-1}\big) - \frac{s-1}{r-1}\log p \le 0$;
\item[(v)] $|\alpha_j(s-1,r) - (j-1)\frac{s-1}{r-1}| < \frac{3}{p\log p}$ for $1\le j \le r$.
\end{itemize}

\end{proposition}

{\sc Proof. }
(i)  It follows from the definition of $\tilde{\mu}:=\tilde{\mu}_p(s-1,r)$ in Corollary \ref{c3} that it satisfies the identity 
$ p^{s-1}\tilde{\mu}^{r-1} = (1+\tilde{\mu})^{r-3}$, clearly implying $0<\tilde{\mu} <1$. For $r=3$, this means that $\tilde{\mu}=\delta$.
For $r\ge 4$, we obtain by virtue of binomial expansion
$$ \tilde{\mu}  = \delta (1+\tilde{\mu})^{\frac{r-3}{r-1}}= \delta + \delta \sum_{k=1}^{\infty} \genfrac{(}{)}{0pt}{}{\frac{r-3}{r-1}}{k} \tilde{\mu}^k,$$ 
where the infinite series has alternating decreasing terms. Hence
$$ 0 < \tilde{\mu} - \delta < \delta\cdot \frac{r-3}{r-1}\cdot \tilde{\mu} < \delta\tilde{\mu},$$
and consequently $\tilde{\mu}< \delta/(1-\delta)$, which implies (i).

(ii) Since the real function $x\mapsto x/(1+x)$ is strictly increasing for $x>0$, we obtain by (i) 
$$ \delta - \frac{\delta^2}{1+\delta} = \frac{\delta}{1+\delta} \le \frac{\tilde{\mu}}{1+\tilde{\mu}} < \frac{\delta + \frac{\delta^2}{1-\delta}}{1+\delta + \frac{\delta^2}{1-\delta}}=\delta.$$
The definition in Corollary \ref{c3} yields that $\tilde{\nu}:=\tilde{\nu}_p(s-1,r) =\tilde{\mu}/(1+\tilde{\mu})$, which proves our claim.

(iii) By (ii), we have $0< \tilde{\nu}< \frac{1}{p}$.
Taking logarithms in the identity $p^{s-1}\tilde{\nu}^{r-1} = (1-\tilde{\nu})^2$ (cf. Cor. \ref{c3}), we obtain
\begin{equation}
\log \frac{1}{\tilde{\nu}} = \frac{s-1}{r-1} \log p - \frac{2}{r-1} \log(1-\tilde{\nu}). \label{62}
\end{equation}
Since $\tilde{\nu} < 1/p$, the Taylor expansion of $\log(1-\tilde{\nu})$ yields
\begin{eqnarray*}
0 < -\log(1-\tilde{\nu}) &=& \sum_{k=1}^{\infty} \frac{\tilde{\nu}^k}{k}
  < \tilde{\nu}+ \frac{\tilde{\nu}^2}{2} \sum_{k=0}^{\infty} \tilde{\nu}^k = \tilde{\nu} + \frac{\tilde{\nu}^2}{2(1-\tilde{\nu})}<\frac{3}{2}\tilde{\nu} < \frac{3}{2p}.
\end{eqnarray*}   
Inserting this into (\ref{62}), we get
$$  0 < \log \frac{1}{\tilde{\nu}} - \frac{s-1}{r-1}\log p < \frac{3}{(r-1)p}.$$  

(iv) For the numbers $\alpha_j:=\alpha_j(s-1,r)$, as defined for $1\le j \le r$ in Corollary \ref{c3}, we have
\begin{equation} 
\alpha_{j+1}-\alpha_j = \left\{\begin{array}{cl}
                           \frac{\log(1/\tilde{\mu})}{\log p} & \mbox{for $j=1$ and $j=r-1$,}\\[.1in]
                           \frac{\log(1/\tilde{\nu})}{\log p} & \mbox{for $2\le j \le r-2$.}
                           \end{array} \right.    \label{diffaj}
\end{equation}                                                      
This is trivial except for $j=r-1$, where it follows from the identities $\tilde{\mu}=\tilde{\nu}/(1-\tilde{\nu})$ and
$p^{s-1}\tilde{\nu}^{r-1}= (1-\tilde{\nu})^2$. Therefore,
$$s-1 = \alpha_r- \alpha_1 = \sum_{j=1}^{r-1} (\alpha_{j+1}-\alpha_j) = \frac{2 \log(1/\tilde{\mu})}{\log p} + \frac{(r-3)\log(1/\tilde{\nu})}{\log p},$$
hence
$$ \log \frac{1}{\tilde{\mu}} - \frac{s-1}{r-1}\log p = \frac{r-3}{2} \left(\frac{s-1}{r-1} \log p -  \log \frac{1}{\tilde{\nu}}\right).$$
Combining this with the bounds found in (iii) completes the argument.

(v) By the definition of the $\alpha_j$, we obtain for $\mbox{$2\le j \le r-1$}$
$$\alpha_j(s-1,r)- (j-1)\frac{s-1}{r-1} =
\left(\frac{\log\frac{1}{\tilde{\mu}}}{\log p} - \frac{s-1}{r-1}\right) + (j-2) \left(\frac{\log\frac{1}{\tilde{\nu}}}{\log p} -\frac{s-1}{r-1}\right).$$
From (iii) and (iv) it follows that
$$ -\,\frac{3}{2p \log p}  < \alpha_j(s-1,r)- (j-1)\frac{s-1}{r-1} <    (j-2) \frac{3}{(r-1)p \log p} < \frac{3}{p \log p},$$
which implies (v) in these cases. Since $\alpha_1=0$ and $\alpha_r=s-1$, the inequality is valid for all $j$.
\begin{flushright}$\Box$\end{flushright}

Proposition \ref{p51}(v) reveals that picking the $\alpha_j$ for $j=1,\ldots,r$  well-spaced in the interval $[0,s-1]$, i.e.
$\alpha_j := (j-1)\frac{s-1}{r-1}$ (see concluding remarks of section 2), is close to best possible. Since it is our
task to find \textit{integral} $a_j$ in optimal position,  it suggests itself to choose the $a_j$ as nearest integers to
the $\alpha_j(s-1,r)$ (as defined in Corollary \ref{c3}) or to the numbers $(j-1)\frac{s-1}{r-1}$, which does not make much of
a difference by Proposition \ref{p51}(v). Anyway, we shall take $a_j = \| \alpha_j \|$ ($1\le j \le r$) with the nearest integer
function $\|\cdot\|$ and have to accept variations between $\alpha_j$ and $a_j$ in the range from $-\frac{1}{2}$ to $\frac{1}{2}$. 
Finally, we shall try to maximize the energy with respect to $r$. \\

We now show that our \textit{integral minimum} $m_p(s-1,r)$, as defined in (\ref{intmin}), can be bounded by the \textit{real minimum} $\tilde{m}_p(s-1,r)$,
introduced in (\ref{realmin}). In general, that is to say in worst cases, we cannot expect to lose less than a factor $p$ between the two minima, taking into account that the shifts from real numbers
$\alpha_j$ to integral parameters $a_j$, varying over an interval of length up to $1$, have to be executed in $h_p$, i.e. in the exponent of $p$.

\begin{proposition}{\label{p52}} 
Let $3\le r\le s$ be given integers, and let $p$ be a prime.
\begin{itemize}
\item[(i)] Let $(\alpha_1,\ldots,\alpha_r) \in \mathbb{R}^r$ be the unique minimizer of $h_p$ determined in Cor. \ref{c3}. 
Then $(a_1,\ldots,a_r) \in \mathbb{Z}^r$ with the nearest integers $a_j:=\|\alpha_j\|$ ($1\le j\le r$) has the property
$$  h_p(a_1,\ldots,a_r) \le \left\{
                            \begin{array}{ll}
                            4\cdot \tilde{m}_2(s-1,r) &  \mbox{ for $p=2$,}  \\
                            p\cdot \tilde{m}_p(s-1,r) &  \mbox{ for $p\ge 3$.} 
                            \end{array}
                            \right.
$$                            
\item[(ii)] We have\quad $\tilde{m}_2(s-1,r) \le m_2(s-1,r) \le  4\cdot \tilde{m}_2(s-1,r).$
\item[(iii)] For any prime $p\ge 3$,
we have 
$$\tilde{m}_p(s-1,r) \le m_p(s-1,r) \le p\cdot \tilde{m}_p(s-1,r).$$
\end{itemize}
\end{proposition}

{\sc Proof. }
The lower bounds in (ii) and (iii) are trivial, and the upper bounds follow immediately from the definition of 
$m_p(s-1,r)$ in (\ref{intmin}). Hence it suffices to prove (i). By Corollary \ref{c3}, we have $\alpha_j = \alpha_j(s-1,r)$ 
for $1\le j \le r$ and
$$ h_p(\alpha_1,\ldots,\alpha_r) = \tilde{m}_p(s-1,r).$$
Since $\tilde{\nu}:=\tilde{\nu}_p(s-1,r) < 1/p$ by Proposition \ref{p51}(ii), we have $\log (1/ \tilde{\nu}) \ge \log p$. 
From Proposition \ref{p51}(i) it follows for $p\ge 3$ that $\tilde{\mu}:=  \tilde{\mu}_p(s-1,r) < 1/\sqrt{p}$,
hence $\log (1/ \tilde{\mu}) > \frac{1}{2}\log p$. 
By use of (\ref{diffaj}), these inequalities imply that $\alpha_2> \frac{1}{2}$, $\alpha_{r-1}<s-1-\frac{1}{2}$ and $\alpha_{j+1} \ge \alpha_j +1$ for $2\le j \le r-2$. 
Moreover, $\alpha_1=0$ and $\alpha_r = s-1$. 
Therefore, the nearest integers $a_j := \|\alpha_j\|$, $1\le j \le r$,  are pairwise distinct, forming a 
strictly increasing sequence. We have $a_j = \alpha_j + \delta_j$ for suitable real numbers $\delta_j$ satisfying $|\delta_j|\le 1/2$ ($1\le j \le r$)
and obtain
\begin{eqnarray*}
h_p(a_1,\ldots,a_r) &=& \sum_{k=1}^{r-1} \sum_{i=k+1}^r \frac{1}{p^{a_i-a_k}} \\
&=& \sum_{k=1}^{r-1} \sum_{i=k+1}^r \frac{1}{p^{(\alpha_i-\delta_i)-(\alpha_k-\delta_k)}}
= \sum_{k=1}^{r-1} \sum_{i=k+1}^r \frac{1}{p^{\delta_k-\delta_i}}\frac{1}{p^{\alpha_i-\alpha_k}}  \\
&\le& p  \sum_{k=1}^{r-1} \sum_{i=k+1}^r \frac{1}{p^{\alpha_i-\alpha_k}}  = p\cdot \tilde{m}_p(s-1,r).
\end{eqnarray*}     
For the prime $p=2$ the above proof has to be modified, since possibly $\alpha_2<\frac{1}{2}$. In this case we choose $a_2=1$ and
$a_{r-1} =s-2$. As before, $\alpha_{j+1} \ge \alpha_j +1$ for $2\le j \le r-2$. Hence we can select each $a_j$, $3\le j \le r-2$, as one of the neighbouring
integers of $\alpha_j$ in such a way that $a_1<a_2<\ldots <a_r$. It follows in this case that $a_j = \alpha_j + \delta_j$ for suitable real numbers $\delta_j$ satisfying 
$|\delta_j|\le 1$ ($1\le j \le r$). Consequently
\begin{equation*}
h_2(a_1,\ldots,a_r) = \sum_{k=1}^{r-1} \sum_{i=k+1}^r \frac{1}{2^{\delta_k-\delta_i}}\frac{1}{2^{\alpha_i-\alpha_k}}
\le 2^2  \sum_{k=1}^{r-1} \sum_{i=k+1}^r \frac{1}{2^{\alpha_i-\alpha_k}}  = 4\cdot \tilde{m}_2(s-1,r).
\end{equation*}
\vspace{-.6truecm}
\begin{flushright}$\Box$\end{flushright}

\setlength{\parindent}{0pt}
{\bf Remark.} The factor $p$ (or $4$ in case $p=2$, respectively) we lose between the \textit{real minimum} $\tilde{m}_2(s-1,r)$
and the \textit{integral minimum} $m_p(s-1,r)$ according to (ii) and (iii) reflects the hypothetical worst case scenario where
each $a_j$ differs from $\alpha_j$ by $\frac{1}{2}$. In practice, the factor between the two minima will be substantially smaller
in almost all cases.  
\vspace{.4in}

\setlength{\parindent}{.6truecm}

In Theorem \ref{t3} the maximal energy $\Emax{p^s}$ as well as the
corresponding $p^s$-maximal sets are given for all primes $p$ and each $s\le 4$, and could be determined quite
easily for other small values of $s$ by (\ref{ft3}), i.e. Theorem 2.1 in \cite{SA2}. 
The inequality
$$   \Emaxstrich{p^s}:= \frac{\Emax{p^s}}{2(p-1)p^{s-1}} \le s$$
is an immediate consequence of Theorem \ref{t4}.
The following result shows that this trivial upper bound lies close to the true value of $\Emaxstrich{p^s}$.

More precisely, part (ii) of the following Theorem \ref{t5} provides the explicit construction of a divisor set ${\cal D}_0$ such that
the energy of the graph $\Icg{{\cal D}_0}{p^s}$ 
falls short of the maximal energy $\Emax{p^s}$ among all integral circulant graphs of order $p^s$
essentially by a factor less than $2$. 
The remark preceding Proposition \ref{p52} explains why we cannot expect to find a more precise lower bound in general.
However, the reader should be aware of the fact that we lose a much smaller factor than $2$ between upper and lower 
bound for $\Emaxstrich{p^s}$ in most cases (cf. the remark following Prop. \ref{p52}). We shall comment on this
at the end of the section.

Bound by the tradition of number theory, \textit{log} will denote the natural logarithm.

\begin{theorem}{\label{t5}} Let $p$ be a prime and let $s$ be a positive integer.
\begin{itemize}
\item[(i)] We have
\begin{equation}
\underline{C}\cdot(s-1)\left(1- \frac{\log\log p}{\log p}\right) \le  \Emaxstrich{p^s} \le 
\overline{C}\cdot(s-1)\left(1- \frac{\log\log p}{\log p}\right)  +1,   \label{Theor5}
\end{equation}
where $\overline{C}=1$ for all $p\ge 3$ and $\underline{C}=\frac{1}{2}$ for all $p\ge 17$ as well as for $3\le p \le 13$ in case $s\le 6$. Only for small values of $p$, we have exceptional constants $\overline{C}=\overline{C}(p)$ and $\underline{C}=\underline{C}(p)$, namely $\overline{C}(2)=0.328$, $\underline{C}(2)=0.118$, and in case $s\ge 7$
$$
\underline{C}(p) = \left\{ \begin{array}{cl}
                            0.030   & \mbox{ if $p=3$, } \\
                            0.233   & \mbox{ if $p=5$, } \\
                            0.337   & \mbox{ if $p=7$, } \\
                            0.442  & \mbox{ if $p=11$, } \\
                            0.473  & \mbox{ if $p=13$. } 
                            \end{array} \right.
$$
%

\item[(ii)] Let $r_0$ be the integer uniquely determined by
$$ \frac{s-1}{D(p)} \le r_0 < \frac{s-1}{D(p)} +1,$$
where
$$
D(p) := \left\{ \begin{array}{cl}
                            4.09184   & \mbox{ for $p=2$, } \\[.1in]
                            2(1+ \frac{\log\log p}{\log p})  & \mbox{ for $p \ge 3$, }
                            \end{array} \right.
$$
and define ${\cal D}_0 = \{p^{\|\alpha_j(s-1,r_0)\|}: \; j=1,\ldots, r_0\}$. 
For $p=2$, $s\ge 11$ and for $p\ge 3$, $s\ge 7$, the energy of the graph $\Icg{{\cal D}_0}{p^s}$ lies in the same interval as the one established for $\Emax{p^s}$ in (\ref{Theor5}).
\end{itemize}
\end{theorem}

{\sc Proof. }
By Theorem \ref{t3} we have for all primes $p$
\begin{equation}
 \Emaxstrich{p^s} = \left\{
            \begin{array}{cl}
            1    & \mbox{ for $s=1$, } \\
            1+\frac{1}{p}    & \mbox{ for $s=2$, } \\
            2-\frac{1}{p} +\frac{1}{p^2}    & \mbox{ for $s=3$, } \\
            2+\frac{1}{p^3}    & \mbox{ for $s=4$. } \\
            \end{array} \right.       \label{supto4}
\end{equation}            
We leave it to the reader to check that each of these values lies within the respective bounds stated 
in (\ref{Theor5}). We may therefore assume $s\ge 5$ in the sequel.

We shall first prove the upper bound in (\ref{Theor5}). By virtue of (\ref{ft5}), it suffices to show that
\begin{equation}
\frac{\Emaxr{p^s}{r}}{2(p-1)p^{s-1}} \le \overline{C}\cdot(s-1)\left(1 - \frac{\log\log p}{\log p}\right) +1.   \label{Theor5rs}
\end{equation}
is satisfied  for all $1\leq r\leq s$.
We distinguish three cases. \\
\underline{Case U1}: $r=1$. \newline
By Corollary 2.1(i) in \cite{SA2}, we have $\Emaxr{p^s}{1} = 2(p-1)p^{s-1}$ for all $p$, which implies (\ref{Theor5rs})
immediately.
\newline
\underline{Case U2}: $r=2$. \newline
It follows from (\ref{ft4}) and Proposition \ref{p1}(i)
that
\begin{equation}
\Emaxr{p^s}{2}  =  2(p-1)p^{s-1} \left(2- (p-1)\frac{1}{p^{s-1}}\right)  \label{spec2}
\end{equation}
for all $p$. Since $s\ge 5$, our upper bound in (\ref{Theor5rs}) is valid in this case.
\newline
\underline{Case U3}: $3\le r\le s$. \newline
By Theorem \ref{t4} and Proposition \ref{p51}(i), we have
\begin{equation}
\Emaxr{p^s}{r} \le 2(p-1) p^{s-1} \left( r-\frac{(p-1)(r-1)}{p^{\frac{s-1}{r-1}}} \right)   \label{upbound}
\end{equation}
for all $p$. Therefore, we study for fixed $p$ and $s$ the real function
$$ g(x):= x - \frac{(p-1)(x-1)}{p^{\frac{s-1}{x-1}}} $$
on the interval $3\le x \le s$ with boundary values
\begin{equation}
g(3)=  3-\frac{2(p-1)}{p^{\frac{s-1}{2}}} \quad \mbox{ and } \quad g(s)= \frac{s-1}{p} +1\, . \label{bounds}
\end{equation}

For a maximum of $g$ at $x_0$, say, with $3< x_0 <s$ the derivative
$$ g'(x_0) = 1 - \frac{p-1}{p^{\frac{s-1}{x_0-1}}} \left(1 + \frac{(s-1) \log p}{x_0-1} \right)$$
vanishes necessarily. Substituting $y:= \frac{s-1}{x-1}$, hence $y\ge 1$, we obtain the condition 
\begin{equation}
1 + y_0 \log p = \frac{p^{y_0}}{p-1}  \label{derivzero}
\end{equation}
for $y_0:=\frac{s-1}{x_0-1}$. Since $x_0 = 1+ (s-1)/y_0$, we conclude for $3 \le x \le s$
\begin{equation} 
g(x) \le g(x_0) = \frac{s-1}{y_0}\left( 1 - \frac{p-1}{p^{y_0}}\right) +1.   \label{gmax}
\end{equation}
\newline
\underline{Case U3.1}: $p=2$. \newline
For $p=2$, equation (\ref{derivzero}) has no solution $y_0\ge 1$, i.e. in that case the maximum
of $g$ is attained at $y_0=1$, that is for $x_0=s\ge 5$, which follows by comparison of the boundary values
in (\ref{bounds}).  Therefore we have in case $p=2$ 
$$ g(x) \le g(s) = s - \frac{s-1}{2} = \frac{s-1}{2} +1.$$
By (\ref{upbound}), this immediately implies (\ref{Theor5rs}).
\newline
\underline{Case U3.2}: $p\ge 3$. \newline
Now (\ref{derivzero}) has a unique solution $y_0$ in the interval $1\le y_0 <2$, corresponding
to the unique maximum of $g$. By a few steps of Newton interpolation we obtain for instance that
$y_0\approx 1.527$ for $p=3$ and $y_0\approx 1.673$ for $p=5$.
We shall verify that
\begin{equation}
g(x) \le (s-1)\left(1 - \frac{\log\log p}{\log p} 
\right)   +1 \label{upin}
\end{equation}
on the interval $3 \le x \le s$. This follows easily for $p=3$ and $p=5$ by inserting the respective values of $y_0$
given above into (\ref{gmax}). For each other fixed prime $p\ge 7$, we define the real function 
$$ w(y) =w_p(y) := \frac{p^y}{p-1}- y \log p -1$$
for all $y\ge 1$. By (\ref{derivzero}) we know that $y_0\ge 1$ satisfies $w(y_0)=0$.
Since the derivative 
$$w'(y) = \left(\frac{p^y}{p-1} -1\right) \log p$$
is positive for $y\ge 1$, the function $w(y)$ is strictly increasing. 
For 
$$ y_p := \frac{\log p}{\log p - \log\log p} - \frac{1}{\log p},$$
which is greater than $1$ for $p\ge 7$, we have
\begin{eqnarray*}
w(y_p) &=& \frac{1}{p-1} p^{\frac{\log p}{\log p - \log\log p}}\cdot e^{-1} - \frac{(\log p)^2}{\log p - \log\log p} \\
       &=& \frac{1}{e(p-1)} p^{1 + \frac{\log\log p}{\log p - \log\log p}} - \frac{(\log p)^2}{\log p - \log\log p} \\
       &<& \frac{p}{e(p-1)} p^{\frac{\log\log p}{\log p - \log\log p}} - \log p < 0\, ,
\end{eqnarray*}
where the final inequality is shown to hold for all primes $p\ge 7$ by simply taking logarithms in
$$ \frac{p}{e(p-1)} p^{\frac{\log\log p}{\log p - \log\log p}} < \log p\, .$$ 
Since $w(y)$ is strictly increasing on $y\ge 1$ and $w(y_p)<0$, but $w(y_0)=0$, it follows that $y_p<y_0$.
By definition of $y_p$, this inequality implies
$$ \log p < \left(1-\frac{\log\log p}{\log p}\right)(1+y_0\log p).$$
Multiplying with $y_0$ and dividing by $(1+y_0\log p)$, we obtain by (\ref{derivzero})
$$ 1 - \frac{p-1}{p^{y_0}} = 1 - \frac{1}{1+y_0\log p} < y_0\left(1-\frac{\log\log p}{\log p}\right).$$
Inserting this into (\ref{gmax}), we have verified (\ref{upin}). By (\ref{upbound}), the proof of
(\ref{Theor5rs}) is complete. Hence the upper bound in (\ref{Theor5}) holds in all cases.

Now we turn our attention to the lower bound for $ \Emaxstrich{p^s}$ an distinguish several cases
and subcases.\\
\underline{Case L1}: $p\ge 3$.\newline
\underline{Case L1.1}: $s\le 4$.\newline
The lower bound in (\ref{Theor5}) has already been verified for $s\le 4$ in (\ref{supto4}).
\newline
\underline{Case L1.2}: $s=5$.\newline
Picking $r=2$, we use (\ref{spec2}) once more and obtain
$$  \Emaxstrich{p^5} \ge \frac{\Emaxr{p^5}{2}}{2(p-1)p^4} = 2- (p-1)\frac{1}{p^{4}} \ge 
2 \left( 1-\frac{\log\log p}{\log p} \right)$$
for all $p\ge 3$, which proves the lower bound of (\ref{Theor5}) in this case. 
\newline
\underline{Case L1.3}: $s=6$.\newline
It follows from (\ref{ft4}) and Proposition \ref{p1}(ii) that
\begin{equation*}
\Emaxstrich{p^6} \ge \frac{\Emaxr{p^6}{3}}{2(p-1)p^5} = 3- (p-1)
\left(\frac{1}{p^2} + \frac{1}{p^5} + \frac{1}{p^3}\right) \ge \frac{5}{2}\left( 1-\frac{\log\log p}{\log p} \right)
\end{equation*}
for all $p\ge 3$, and again the lower bound in (\ref{Theor5}) is confirmed.
\newline
\underline{Case L1.4}: $s\ge 7$.\newline
We choose the integer $r_0$ according to the inequality
\begin{equation}
\frac{s-1}{2L_p} \le r_0 < \frac{s-1}{2L_p} +1,  \label{rnull1}
\end{equation}
where 
$$L_p := 1+ \frac{\log\log p}{\log p}.$$
A simple calculation reveals that for $s\ge 7$ and all primes $p\ge 3$ the expression on the left-hand
side of (\ref{rnull1}) is always greater than $2$. Consequently, $3 \le r_0 \le s-1$.
By Corollary \ref{c3}, we have
\begin{equation}
\tilde{m}_p(s-1,r_0) = (r_0-1 + \tilde{\mu}_p(s-1,r_0))\cdot \tilde{\mu}_p(s-1,r_0),   \label{uppm1}
\end{equation}
and from Proposition \ref{p51}(i) and the definition of $r_0$, we get $\tilde{\mu}_p(s-1,r_0)<\delta +\delta^2/(1-\delta)$ for $\delta = p^{-(s-1)/(r_0-1)}<p^{-2L_p}$. Hence
$$ \tilde{\mu}_p(s-1,r_0)< \frac{1}{p^{2L_p}}\left(1+\frac{1}{p^{2L_p}-1}\right).$$
By (\ref{uppm1}) and (\ref{rnull1}) we now have
\begin{eqnarray*}
\tilde{m}_p(s-1,r_0) &<& \left( \frac{s-1}{2L_p}+ \frac{1}{p^{2L_p}}\left(1+\frac{1}{p^{2L_p}-1}\right) \right)
    \cdot \frac{1}{p^{2L_p}}\left(1+\frac{1}{p^{2L_p}-1}\right)  \\
    &<& \left( \frac{s-1}{2p^{2L_p}L_p}+ \frac{1}{p^{4L_p}}\right)\cdot \frac{p}{p-1}\, ,
\end{eqnarray*}
because for $p\ge 3$
\begin{eqnarray*}
1+\frac{1}{p^{2L_p}-1} &<& \left(1+\frac{1}{p^{2L_p}-1}\right)^2 = \left(1+\frac{1}{p^2(\log p)^2-1}\right)^2 \\
&<& \left(1+\frac{1}{(p-1)(p+1)}\right)^2 < \frac{p}{p-1}\,.
\end{eqnarray*}
Since $p\ge 3$, we obtain by Proposition \ref{p52}(iii) that
\begin{eqnarray*}
m_p(s-1,r_0) &<& p\cdot \tilde{m}_p(s-1,r_0) \\
&<& \left( \frac{s-1}{2p^{2L_p}L_p}+ \frac{1}{p^{4L_p}}\right)\cdot \frac{p^2}{p-1} \\
&=& \left(\frac{s-1}{2(\log p)^2L_p} +
\frac{1}{p^2(\log p)^4 L_p}\right) \cdot \frac{1}{p-1}\,. 
\end{eqnarray*}
According to (\ref{ft4}) and (\ref{rnull1}), it follows that
\begin{eqnarray}
\frac{\Emaxr{p^s}{r_0}}{2(p-1)p^{s-1}} &=& r_0 -(p-1)\,m_p(s-1,r_0) \nonumber \\
&>& \frac{s-1}{2L_p} - \left(\frac{s-1}{2(\log p)^2L_p} +\frac{1}{p^2(\log p)^4}\right) \label{lowbound1}\\ 
&>&  \frac{s-1}{2L_p}\left(1-\frac{1}{(\log p)^2}\right) - \frac{1}{p^2(\log p)^4}\, . 
\nonumber
\end{eqnarray}
It is easy to check that for all primes $p\ge 17$
$$ (\log\log p)^2 > 1+ \frac{L_p}{3(p\log p)^2},$$
and that for $3\le p \le 13$
$$ (\log\log p)^2 > 1+ \left(\frac{1}{3(p\log p)^2}- c_p(\log p)^2\right) L_p $$
with constants $c_p$ defined by the following table:
\begin{center}
\begin{tabular}{|c|c|} \hline
$p$ & $c_p$ \\ \hline
& \\[-.18in]
{\rm 3} & {\rm 0.859} \\
{\rm 5} & {\rm 0.375} \\
{\rm 7} & {\rm 0.214}\\
{\rm 11}& {\rm 0.073}\\
{\rm 13}& {\rm 0.033}\\ \hline
\end{tabular}
\end{center}
Setting $c_p:=0$ for all primes $p\ge 17$, the fact that $s\ge 7$ implies for all $p\ge 3$
$$  (\log\log p)^2 + c_p\log p\, L_p > 1 + \frac{L_p}{3(p\log p)^2} 
\ge 1 +\frac{2L_p}{(s-1)(p\log p)^2}.$$
Dividing by $(\log p)^2$, adding $1$ on both sides and rearranging terms yields
\begin{eqnarray*}
1- \frac{1}{(\log p)^2} - \frac{2L_p}{(s-1)p^2(\log p)^4}
&>& 1- \left(\frac{\log\log p}{\log p}\right)^2 - \frac{c_pL_p}{\log p} \\
&=& \left(1- \frac{\log\log p}{\log p}\right)L_p - 
\frac{c_pL_p}{\log p}.
\end{eqnarray*}
Dividing by $L_p$ and multiplying with $(s-1)/2$ implies
$$  \frac{s-1}{2L_p}\left(1- \frac{1}{(\log p)^2}\right) -\frac{2}{(s-1)p^2(\log p)^4}
> \frac{s-1}{2}\left(1- \frac{\log\log p + c_p}{\log p}\right).$$
Inserting this inequality into (\ref{lowbound1}) yields
$$ \frac{\Emaxr{p^s}{r_0}}{2(p-1)p^{s-1}} > \frac{s-1}{2} \left(1-\frac{\log\log p +c_p}{\log p}\right),$$
which completes the proof of the lower bound in (\ref{Theor5}) for $p\ge 3$. 
At the same time, our construction combined with Proposition \ref{p52}(i) reveals the truth of statement (ii) for primes $p\ge 3$ and $s\ge 7$.\\
\newline
\underline{Case L2}: $p=2$. \newline
\underline{Case L2.1}: $s\le 4$.\newline
The lower bound in (\ref{Theor5}) follows from (\ref{supto4}).
\newline
\underline{Case L2.2}: $5\le s\le 10$.\newline
Picking $r=3$, it follows from (\ref{ft4}) and Proposition \ref{p1}(ii) that
\begin{equation*}
\Emaxstrich{2^s} \ge \frac{\Emaxr{2^s}{3}}{2^s} = 
3- \left(\frac{1}{2^{[\frac{s-1}{2}]}} + \frac{1}{2^{s-1}} + \frac{1}{2^{s-1-[\frac{s-1}{2}]}}\right) .
\end{equation*}
It is easy to check that the last term becomes minimal for $s=5$. Hence
$$ \Emaxstrich{2^s} \ge 3 - \frac{9}{16} \ge  \frac{2}{11} (s-1)$$
for all $s$ in the given range. This confirms the lower bound in (\ref{Theor5}) for these values of $s$.
\newline
\underline{Case L2.3}: $s\ge 11$.\newline
Let $c_1>1$ be the unique real number satisfying $c_1^2-6c_1+5 =4c_1\log c_1$, i.e. $17.0517 < c_1 < 17.0518$, and
let $c_2:= \frac{\log c_1}{\log 2}$, thus $4.09184 < c_2 < 4.09186$.
We choose the integer $r_2$ according to the inequality
\begin{equation}
\frac{s-1}{c_2} \le r_2 < \frac{s-1}{c_2} +1.  \label{rzwei1}
\end{equation}
Apparently, the expression on the left-hand side of (\ref{rzwei1}) is always greater than $2$ for $s\ge 11$. Consequently, $3 \le r_2 \le s-1$.
By Corollary \ref{c3}, we have
\begin{equation}
\tilde{m}_2(s-1,r_2) = (r_2-1 + \tilde{\mu}_p(s-1,r_2))\cdot \tilde{\mu}_p(s-1,r_2),   \label{uppm2-1}
\end{equation}
and from Proposition \ref{p51}(i) and the definition of $r_2$, we get 
$$ \tilde{\mu}_p(s-1,r_2)<  \frac{1}{2^{c_2}}+ \frac{1}{2^{2c_2}-2^{c_2}} = \frac{1}{2^{c_2}}\left(1+\frac{1}{2^{c_2}-1}\right).$$
By (\ref{uppm2-1}) and (\ref{rzwei1}) we obtain
\begin{equation*}
\tilde{m}_2(s-1,r_2) < \left( \frac{s-1}{c_2}+ \frac{1}{2^{c_2}}\left(1+\frac{1}{2^{c_2}-1}\right) \right)
    \cdot \frac{1}{2^{c_2}}\left(1+\frac{1}{2^{c_2}-1}\right).
\end{equation*}
Proposition \ref{p52}(ii) implies that
$  m_2(s-1,r_2) < 4\cdot \tilde{m}_2(s-1,r_2),$
and with (\ref{ft4}) and (\ref{rzwei1}) we get
\begin{eqnarray*}
\Emaxstrich{2^s} &\ge& \frac{\Emaxr{2^s}{r_2}}{2^s} = r_2 - m_2(s-1,r_2)  \\
&>& \frac{s-1}{c_2} - 4 \left(\frac{s-1}{c_2}+ \frac{1}{2^{c_2}}\left(1+\frac{1}{2^{c_2}-1}\right) \right)
    \cdot \frac{1}{2^{c_2}}\left(1+\frac{1}{2^{c_2}-1}\right)    \\ 
&=&  \frac{s-1}{c_2}\left(1-\frac{4}{2^{c_2}}\left(1+\frac{1}{2^{c_2}-1}\right) \right) - \frac{4}{2^{2c_2}}\left(1+\frac{1}{2^{c_2}-1}\right)^2 \\
&>& \frac{s-1}{5.45} - 0.01553 > \frac{2}{11} (s-1)
\end{eqnarray*}
for $s\ge 11$. This proves the lower bound of (\ref{Theor5}) and completes part (i) of Theorem \ref{t5}. This time our construction combined with Proposition \ref{p52}(i) shows (ii) for $p=2$ and $s\ge 11$.
\vspace{-.2truecm}
\begin{flushright}$\Box$\end{flushright}
\vspace{.2in}

As an example, the following table illustrates the previous theorem for $s=17$ and several small values of $p$.
Note that $r_0$ will eventually become $8$, roughly for $p>10^{10}$).

\begin{center}
\begin{tabular}{cccrrr}
{$n$} & {$r_0$} 
& {${\cal D}_0$} & lower & $\Emaxstrich{n}$ &  upper \\ \hline
  $3^{17}$ & $8$ 
        & $(0, 2, 5, 7, 9, 11, 14, 16)$
        & $0.439$ & $6.652$ & $15.630$ \\ 
  $5^{17}$ & $7$ 
        & $(0, 3, 5, 8, 11, 13, 16)$
        & $2.626$ & $6.547$ & $12.269$ \\ 
  $7^{17}$ & $6$ 
        & $(0, 3, 6, 10, 13, 16)$
        & $3.547$ & $5.927$ & $11.526$ \\ 
  $11^{17}$ & $6$ 
        & $(0, 3, 6, 10, 13, 16)$
        & $4.493$ & $5.969$ & $11.164$ \\ 
  $13^{17}$ & $6$ 
        & $(0, 3, 6, 10, 13, 16)$
        & $4.789$ & $5.978$ & $11.124$ \\ 
  $17^{17}$ & $6$ 
        & $(0, 3, 6, 10, 13, 16)$
        & $5.059$ & $5.987$ & $11.119$ \\ 
  $23^{17}$ & $6$ 
        & $(0, 3, 6, 10, 13, 16)$
        & $5.084$ & $5.993$ & $11.169$ \\ 
\end{tabular}
\end{center}

\medskip

\setlength{\parindent}{0pt}
{\bf Remarks.}
\begin{itemize}
\item[(i)]
The proof of Theorem \ref{t5} shows that we lose a factor $2$ between lower and upper bound for
$\Emaxstrich{p^s}$ in case $p\ge 17$ (and similarly for smaller $p$) mainly due to the fact that we lose a factor
$p$ between $\tilde{m}_p(s-1,r)$ and $m_p(s-1,r)$, which however happens only as an extremely rare worst case event (cf. Proposition \ref{p52} and the preceding and subsequent remarks). A staightforward adaptation of the method introduced in the proof for the lower bound of (\ref{Theor5}) implies the following:

Assume that for some fixed sufficiently large $p$ and $s$ we have $m_p(s-1,r) \le p^{\gamma}\cdot \tilde{m}_p(s-1,r)$ with some positive
$\gamma <1$. Taking 
\begin{equation*}
\frac{s+1}{(1+\gamma)L_p} \le r_0 < \frac{s+1}{(1+\gamma)L_p} +1
\end{equation*}
instead of (\ref{rzwei1}), we obtain
\begin{equation*}
\frac{1}{1+\gamma}(s+1)\left(1-\frac{\log\log p}{\log p}\right) \le  \Emaxstrich{p^s}\le (s+1)\left(1 - \frac{\log\log p}{\log p}\right) +1. 
\end{equation*}
This reveals that, the smaller the difference between integral and real maximum is, the better our bounds are. In the extreme case where $m_p(s-1,r)= \tilde{m}_p(s-1,r)$, i.e.
$\gamma =0$, lower and upper bound differ only by $1$, and the gcd graph with the corresponding divisor set ${\cal D}_0$ has maximal energy, because the energy is an integral 
number, but the two  bounds are not.
\item[(ii)]
It should be noted that even the lower bound of (\ref{Theor5}) already implies hyperenergeticity in most cases. A straightforward 
calculation shows this for e.g.~$p\geq 3$ and $s\geq 5$.
\end{itemize}

\setlength{\parindent}{.6truecm}

%
%

\section{Conclusion and open problems}

Given a fixed prime power $p^s$, we have provided a method to construct a divisor set ${\cal D}_0$ with the
property that $\Ene{{\cal D}_0}{p^s}\ge \frac{1}{2} \Emax{p^s}$. In most cases we expect $\Ene{{\cal D}_0}{p^s}$
to lie much closer to $\Emax{p^s}$ than our worst case inequality guarantees. But since we have used the
``real" maximum for reference it may not be expected to get hold of the ``integral" maximum in general, using an analytic approach.
The convexity properties of the function $h_p$ also suggest that a divisor set  ${\cal D}_{\max}$ with $\Ene{{\cal D}_{\max}}{p^s} = \Emax{p^s}$ can be found in the ``neighborhood" of our ${\cal D}_0$. Given an explicit
integer $p^s$ it should not be too difficult to determine $\Emax{p^s}$ precisely by comparison of a
few candidates for a ${\cal D}_{\max}$ ``near" ${\cal D}_0$.

Let us conclude this section by posing the challenge of finding similarly accessible bounds on
$\Emax{n}$ for integers $n$ which have different prime factors. Even for $n=p_1^{s_1}p_2^{s_2}$ with primes
$p_1\neq p_2$ and arbitrary divisor sets a closed formula for the energy of the corresponding integral circulant graphs would
be most desirable. Of course, this should then be the basis for analyzing these graphs for hyperenergeticity.

%
%

\end{document}